\input amstex
\loadbold
\loadeurm
\loadeusm
\openup.8\jot \magnification=1200

\def\L{\Lambda}
\def\O{\Omega}

\def\Z{\Bbb Z}

\def\det{{\text{det}}}

\def\H{{\Bbb H}}

\def\C{{\Bbb C}}

\def\G{{\Cal G}}
\def\Y{{\Cal Y}}
\def\P{{\Cal P}}

\documentstyle{amsppt}
\pageheight{7.7in}
\vcorrection{-0.05in}
\topmatter
\pageno 1
\title Explicit matrices for Hecke operators
on Siegel modular forms\endtitle
\author Lynne H. Walling\endauthor
\address L.H. Walling, Department of Mathematics,
University of Bristol, Bristol BS8 1TW, England \endaddress
\email l.walling\@bristol.ac.uk \endemail

\abstract We present an explicit set of matrices giving the action of
the Hecke operators $T(p), T_j(p^2)$ on Siegel modular forms.

\endabstract

\endtopmatter
\document

\head{Introduction}\endhead
It is well-known that the space of elliptic modular forms of weight
$k$ has a basis of simultaneous eigenforms for the Hecke operators,
and the Fourier coefficients of an eigenform (and hence the eigenform)
are completely determined by its eigenvalues and first Fourier
coefficient.  In the theory of Siegel modular forms, the role of the
Hecke operators is not yet completely understood, thus there are many
avenues open for conjecture and exploration, including computational
exploration.  The purpose of this note is to present an explicit set
of matrices giving the action of the Hecke operators on Siegel modular
forms, with the goal of facilitating
 computational exploration.  (This construction also yields an
explicit set of matrices giving the action of Hecke operators on
Jacobi modular forms; we remark on this further at the end of this
note.)
\smallskip

\head{Definitions and results}\endhead

For $F$ a Siegel modular form of degree $n$ and $p$ a prime, we define
the Hecke operator $T(p)$ by
$$F|T(p)=p^{n(k-n-1)/2}\sum_{M}F|\pmatrix {1\over p}I_n\\&I_n\endpmatrix M$$
where $M$ runs over a complete set of coset representatives for
$(\Gamma'\cap\Gamma)\backslash\Gamma$, with $\Gamma=Sp_n(\Z)$ and
$$\Gamma'=\pmatrix pI_n\\&I_n\endpmatrix\Gamma
\pmatrix {1\over p}I_n\\&I_n\endpmatrix.$$
Similarly, for $1\le j\le n$, we define $T_j(p^2)$ by
$$F|T_j(p^2)=
\sum_M F|\pmatrix {1\over
p}I_j\\&I_{n-j}\\&&pI_j\\&&&I_{n-j}\endpmatrix M$$
where $M$ runs over a complete set of coset representatives for
$(\Gamma_j'\cap\Gamma)\backslash\Gamma$; here 
$$\Gamma_j'=\pmatrix pI_j\\&I_{n-j}\\&&{1\over p}I_j\\&&&I_{n-j}\endpmatrix\Gamma
\pmatrix {1\over p}I_j\\&I_{n-j}\\&&pI_j\\&&&I_{n-j}\endpmatrix.$$

In  [3], we determine the action of Hecke operators on Fourier
coefficients of a Siegel modular form by first describing
 a complete set of coset
representatives for the Hecke operators.  
We index the cosets using  lattices;
 the coset representatives are then explicitly described
except for a choice of $G\in GL_n(\Z)$ associated to each lattice.
There are infinitely many possible choices for each $G$; in this note we
make an explicit choice for each $G$.  

 We first construct matrices for $T_j(p^2)$ (and for the averaged
operators $\widetilde T_j(p^2)$ introduced in [3]); then we do the same
for $T(p)$.

For $T_j(p^2)$, we construct these $G$ as follows.
For (nonnegative) integers $r_0,r_2$ with $r_0+r_2\le j$ and
$r_1=j-r_0-r_2$,
 we call $\P$ a
partition of type $(r_0,r_2)$ for $(n,j)$
 if $\P$ is an ordered partition 
$$(\{d_1,\ldots,d_{r_0}\}, \{b_1,\ldots,b_{r_1}\}, \{a_1,\ldots,a_{r_2}\}, 
\{c_1,\cdots,c_{n-j}\})$$
 of $\{1,2,\ldots,n\}$.
(Note that if some $r_i=0$ or $n-j=0$, a set in the partition could be
empty.)  Given a partition $\P$ of type $(r_0,r_2)$, we let
$\G_{\P}\subseteq GL_n(\Z)$
consist of all matrices $G=(G_0,G_1,G_2,G_3)$ constructed as follows.  
$G_0$ is the $n\times r_0$ matrix with $\ell,t$-entry
1 if $\ell=d_t$, and 0 otherwise.
$G_1$ is an $n\times r_1$ matrix with $\ell,t$-entry $\beta_{\ell t}$
where $\beta_{\ell t}=1$ if $\ell=b_t$, $\beta_{\ell t}=0$ if
$\ell<b_t$ or $\ell=a_i$ (some $i$) or $\ell=b_i$ (some $i\not=t$),
and otherwise $\beta_{\ell t}\in\{0,1,\ldots,p-1\}$.  $G_2'$ is an
$n\times r_2$ matrix with $\ell,t$-entry $\alpha_{\ell t}$ where
$\alpha_{\ell t}=1$ if $\ell=a_t$, $\alpha_{\ell t}=0$ if $\ell<a_t$
or $\ell=a_i$ (some $i\not=t$), and otherwise $\alpha_{\ell
t}\in\{0,1,\ldots,p-1\}$.  $G_2''$ is an $n\times r_2$ matrix with
$\ell,t$-entry $\delta_{\ell t}$ where $\delta_{\ell t}=0$ if
$\ell\not=d_i$ (any $i$), and otherwise $\delta_{\ell
t}\in\{0,1,\ldots,p-1\}$.  $G_2=G_2'+pG_2''$.  $G_3$ is an
$n\times(n-j)$ matrix with $\ell,t$-entry $\gamma_{\ell t}$ where
$\gamma_{\ell t}=1$ if $\ell=c_t$, $\gamma_{\ell t}=0$ if $\ell<c_t$
or $\ell=a_i$ or $b_i$ (some $i$) or $\ell=c_i$ (some $i\not=t$), and
otherwise $\gamma_{\ell t}\in\{0,1,\ldots,p-1\}$.

Note that $(G_0,G_1,G_2',G_3)$ is a (column) permutation of an integral lower
triangular matrix with 1's on the diagonal, and thus is an element of
$GL_n(\Z)$.  Also, it is easy to see that there is an elementary
matrix $E$ so that 
$$(G_0,G_1,G_2',G_3)E=(G_0,G_1,G_2'+pG_2'',G_3)=G,$$
and so $G\in GL_n(\Z)$.  (After proving Theorem 1, we describe $G^{-1}$
as a product of four explicit matrices.)

We let $\G_{r_0,r_2}=\cup_{\P} \G_{\P}$ where $\P$ varies over all
partitions of type $(r_0,r_2)$.
We set 
$$D_{r_0,r_2}=\pmatrix
I_{r_0}\\&pI_{r_1}\\&&p^2I_{r_2}\\&&&I_{n-j}\endpmatrix.$$
Also,
we let $\Y_{r_0,r_2}$ be the set of all (integral) matrices of the form
$$\pmatrix Y_0&Y_2&0&Y_3\\p\ ^tY_2&Y_1\\0\\^tY_3\endpmatrix$$
where
$Y_0$ is symmetric, $r_0\times r_0$, with entries varying modulo
$p^2$, $Y_1$ is symmetric, $r_1\times r_1$, with entries varying modulo $p$,
$Y_2$ is $r_0\times r_1$ with entries varying modulo $p$, and $Y_3$ is
$r_0\times (n-j)$ with entries varying modulo $p$.  We let
$\Y'_{r_0,r_2}$ be those matrices in $\Y_{r_0,y_2}$ that satisfy the
additional condition $p\nmid\det Y_1$ (which is trivially satisfied if $r_1=0$).

\proclaim{Theorem 1}  Given a degree $n$ Siegel modular form $F$,
$1\le j\le n$, 
$$F|T_j(p^2)=
\sum_M F|\pmatrix {1\over
p}I_j\\&I_{n-j}\\&&pI_j\\&&&I_{n-j}\endpmatrix M$$
where $$M=\pmatrix D&Y\\&D^{-1}\endpmatrix \pmatrix G^{-1}\\&^tG\endpmatrix$$
varies so that for some $r_0,r_2$ with $r_0+r_2\le j$,
$D=D_{r_0,r_2}$, $Y\in\Y'_{r_0,r_2}$, and $G\in\G_{r_0,r_2}$.
Also, with 
$$\widetilde T_j(p^2)=p^{j(k-n-1)}\sum_{0\le t\le j}\left[{n-t \atop
j-t}\right]_p T_t(p^2)$$ where
$\left[{m\atop r}\right]_p=\prod_{i=0}^{r-1}{p^{m-i}-1\over
p^{r-i}-1}$,
$$F|\widetilde T_j(p^2)=
p^{j(k-n-1)} \sum_M F|\pmatrix {1\over
p}I_j\\&I_{n-j}\\&&pI_j\\&&&I_{n-j}\endpmatrix M$$
where $$M=\pmatrix D&Y\\&D^{-1}\endpmatrix \pmatrix G^{-1}\\&^tG\endpmatrix$$
varies so that for some $r_0,r_2$ with $r_0+r_2\le j$,
$D=D_{r_0,r_2}$, $Y\in\Y_{r_0,r_2}$, and $G\in\G_{r_0,r_2}$.
\endproclaim

\demo{Proof} As mentioned above,
in Proposition 2.1 of [3] we found a complete set of coset
representatives indexed by lattices $(\O,\L_1)$ where, for $\L$ a
fixed lattice of rank $n$, $\O$ varies subject to
$p\L\subseteq\O\subseteq{1\over p}\L$, and $\overline\L_1$ varies over
all codimension $n-j$ subspaces of $\L\cap\O/p(\L+\O)$.  With
$\{\L:\O\}$ denoting the invariant factors of $\O$ in $\L$ and $r_0$
the multiplicity of the invariant factor $p$, $r_2$ the multiplicity
of the invariant factor $1\over p$, the action of the coset representatives
corresponding to $(\O,\L_1)$ is given by the matrices
$$\pmatrix D&Y\\&D^{-1}\endpmatrix \pmatrix G^{-1}\\&^tG\endpmatrix$$
where $D=D_{r_0,r_2}$, $Y$ varies over $\Y'_{r_0,r_2}$, and $G=G(\O,\L_1)$ is any
change of basis matrix so that, relative to a fixed basis
$(x_1,\ldots,x_n)$ for $\L$,
$$\O=\L G D^{-1}\pmatrix pI_j\\&I_{n-j}\endpmatrix,\ 
\L_1=\L G\pmatrix 0_{r_0}\\&I_{r_1}\\&&0\endpmatrix$$
where $r_1=j-r_0-r_2$.  (Thus only those $\O$ occur where $r_0+r_2\le
j$.) So we need to show that such a pair $(\O,\L_1)$ corresponds to a
unique $G\in\G_{r_0,r_2}$.  We do this by constructing all $(\O,\L_1)$
for each pair of parameters $(r_0,r_2)$, simultaneously building a
partition $\P$ and making choices of $\alpha_{\ell t}, \beta_{\ell
t}, \gamma_{\ell t}, \delta_{\ell t}$ so that with $G$ constructed
according to our recipe above, we can take $G$ for $G(\O,\L_1)$.

Notice that when $p\L\subseteq\O\subseteq{1\over p}\L$, the
Invariant Factor Theorem (81:11 of [5])
tells us we have compatible decompositions:
$$\align
\L&=\L_0\oplus\L_1'\oplus\L_2,\\
\O&=p\L_0\oplus\L_1'\oplus{1\over p}\L_2.
\endalign$$
On the other hand, given $\L$, such an $\O$ is determined by $\O'=\L_2+p\L$
and $(p\L_1'\oplus\L_2)+p\O'$.  Also, in $\L\cap\O/p(\L+\O)$, 
$\overline\L_2=0$, so $\L_1$ can be chosen so that in $\O'/p\O'$,
$\overline\L_1\subseteq\overline{p\L'_1}\subseteq\overline{p\L}$.

So to begin our construction of $\O,\L_1$ and $G=G(\O,\L_1)$,
in $\L/p\L$ we choose a dimension $r_2$ subspace $\overline
C'$; let $(\overline v_1',\ldots,\overline v_{r_2}')$ be a basis for
$\overline C'$.  Each $\overline v_t'$ is a linear combination over
$\Z/p\Z$ of the $\overline x_i$; 
 by adjusting the $\overline v_t'$ we can assume
$$\overline v_t'=\overline x_{a_t}
+ \sum_{\ell>a_t} \overline \alpha_{\ell t} \overline x_{\ell}$$
where $a_1,\ldots,a_{r_2}$ and distinct and
 $\overline\alpha_{\ell t}=0$ if $\ell=a_i$ (some $i\not=t$).  Let
$\alpha_{\ell t}\in\{0,1,\ldots,p-1\}$ be a preimage of
$\overline\alpha_{\ell t}$.

Now let $\O'$ be the preimage in $\L$ of $\overline C'$.  In
$\O'/p\O'$ we will construct a dimension $n-r_0$ subspace $\overline C$ so
that $\dim(\overline C\cap\overline{p\L})=n-r_0-r_2$, distiguishing a dimension
$r_1$ subspace $\overline{p\L_1}$ of $\overline C\cap\overline{p\L}$.
We begin by choosing $\overline {p\L_1}$ to be a dimension $r_1$
subspace of $\overline{p\L}$; let
$\overline{pu_1},\ldots,\overline{pu_{r_1}}$ be a basis for
$\overline{p\L_1}$.  Since $\overline{px_{a_i}}=0$ in $\O'/p\O'$,
we can adjust the $\overline{pu_t}$ so that
$$\overline {pu_t}=\overline{px_{b_t}}
+ \sum_{\ell>b_t} \overline\beta_{\ell t}\overline{px_{\ell}}$$
where $b_1,\ldots,b_{r_1}$ are distinct, $b_t\not=a_i$ (any $i$),
and $\overline\beta_{\ell t}=0$ if $\ell=a_i$ (some $i$)
or $\ell=b_i$ (some $i\not=t$).  Let
$\beta_{\ell t}\in\{0,1,\ldots,p-1\}$ be a preimage of
$\overline\beta_{\ell t}$.

Now extend $\overline{p\L_1}$ to a dimension $n-r_0-r_2$ subspace
$\overline{p \L_1'}$ of $\overline{p\L}$ in $\O'/p\O'$.
Extend $(\overline{pu_1},\ldots,\overline{pu_{r_1}})$ to a basis
$$(\overline{pu_1},\ldots,\overline{pu_{r_1}},
\overline{pw_1},\ldots,\overline{pw_{n-j}})$$
for $\overline{p\L_1'}$ so that
$$\overline{pw_t}=\overline{px_{c_t}}+\sum_{\ell>c_t}\overline\gamma_{\ell
t}\overline{px_{\ell}},$$ 
where $c_1,\ldots,c_{n-j}$ are distinct,
$c_t\not=a_i, b_i$ (any $i$), and $\overline\gamma_{\ell t}=0$
if $\ell=a_i$ (some $i$), or $\ell=b_i$ (some $i\not=t$).  Let
$\gamma_{\ell t}\in\{0,1,\ldots,p-1\}$ be a preimage of
$\overline\gamma_{\ell t}$.  

Now we extend $\overline{p\L_1'}$ to a dimension $n-r_0$ space
$\overline C$ so that the dimension of $\overline C\cap\overline{p\L}$
is $n-r_0-r_2=r_1+n-j$, and we extend
$(\overline{pu_1},\ldots,\overline{pw_1},\ldots)$ to a
basis $$(\overline{pu_1},\ldots,\overline{pu_{r_1}},\overline{pw_1},\ldots,\overline{pw_{n-j}},\overline{pv_1},\ldots,\overline{pv_{r_2}})$$
for $\overline C$.  Taking $d_1,\ldots,d_{r_0}$ so that
$$(\{d_1,\ldots,d_{r_0}\},\{b_1,\ldots,b_{r_1}\},\{a_1,\ldots,a_{r_2}\},
\{c_1,\ldots,c_{n-j}\})$$
is a partition of $\{1,\ldots,n\}$, we can take
$$\overline v_t=\overline v_t' 
+ \sum_{m=1}^{r_0}\overline\delta_{m t}\overline{px_{d_m}}$$
for some $\overline\delta_{mt}$;
let $\delta_{mt}\in\{0,1,\ldots,p-1\}$ be a preimage of
$\overline\delta_{mt}$.  

Now let $p\O$ be the preimage in $\O'$ of
$\overline C$.  So with
$$\align
u_t&=x_{b_t}+\sum_{\ell>b_t}\beta_{\ell t}x_{\ell}\ (1\le t\le r_1),\\
v_t&=x_{a_t}+\sum_{\ell>a_t}\alpha_{\ell
t}x_{\ell}+p\sum_{m}\delta_{mt}x_{d_m}\ (1\le t\le r_2),\\
w_t&=x_{c_t}+\sum_{\ell>c_t}\gamma_{\ell t}x_{\ell}\ (1\le t\le n-j),
\endalign$$
the vectors
$$(px_{d_1},\ldots,px_{d_{r_0}},u_1,\ldots,u_{r_1},{1\over
p}v_1,\ldots,
{1\over p}v_{r_2},w_1,\ldots,w_{n-j})$$
form a basis for $\O$, and $(\overline u_1,\ldots,\overline u_{r_1})$
is a basis for $\overline \L_1$ in $\L\cap\O/p(\L+\O)$.
$\square$
\enddemo

\remark{Remark} Given $G\in\G_{\P}$ as above, $G^{-1}=E_1E_2E_3E_4$
where the $E_i$ are $n\times n$ matrices constructed as follows.
$E_1$ has $i,i$-entry 1 ($1\le i\le n$); for $1\le\ell\le r_0$, $E_1$ has 
$\ell,r_0+t$-entry $-\beta_{d_{\ell}t}$ 
($1\le t\le r_1$), $\ell,r_0+r_1+t$-entry
$-\alpha_{d_{\ell}t}-p\delta_{d_{\ell}t}$ ($1\le t \le r_2$),
$\ell,r_0+r_1+r_2+t$-entry $-\gamma_{d_{\ell}t}$ ($1\le t\le n-j$),
and all other entries 0.  So for $1\le t\le n-j$, column
$r_0+r_1+r_2+t$ of $GE_1$ has a 1 in row $c_t$, and zeros elsewhere.
$E_2$ has $i,i$-entry 1 ($1\le i\le n$); for $1\le \ell\le n-j$, $E_2$ has
$r_0+r_1+r_2+\ell,r_0+t$-entry $-\beta_{c_{\ell}t}$ ($1\le t\le r_1$),
$r_0+r_1+r_2+\ell,r_0+r_1+t$-entry $-\alpha_{c_{\ell}t}$ ($1\le t\le
r_2$), and zeros elsewhere.  Thus for $1\le t\le r_1$, column $r_0+t$
of $GE_1E_2$ has a 1 in row $b_t$, and zeros elsewhere.  $E_3$ has
$i,i$-entry 1 ($1\le i\le n$); for $1\le \ell\le r_1$, $1\le t\le
r_2$, $E_3$ has $r_0+\ell,r_0+r_1+t$-entry $-\alpha_{b_{\ell}t}$ and
zeros elsewhere.  Thus $GE_1E_2E_3=\ ^tE_4$ is a permutation matrix;
$E_4$ has 1 as its
$\ell,d_{\ell}$-entry  ($1\le \ell\le r_0$),  1 as its 
$r_0+\ell,b_{\ell}$-entry ($1\le \ell\le r_1$),  1 as its
$r_0+r_1+\ell,a_{\ell}$-entry ($1\le \ell\le r_2$),  1 as its
$r_0+r_1+r_2+\ell,c_{\ell}$-entry ($1\le \ell\le n-j$), and zeros
elsewhere.  So $GE_1E_2E_3E_4=I$.
\endremark

We follow a similar procedure to construct matrices for $T(p)$: For
$0\le r\le n$, we let $\G_r$ be the set of matrices $G$ so that for some
ordered partition $\P=(\{d_1,\ldots,d_r\},\{a_1,\ldots,a_{n-r}\})$ of
$\{1,2,\ldots,n\}$, for $1\le t\le r$,
 column $t$ of $G$ has  1 in row $d_t$ and zeros elsewhere,
 and for $1\le t\le n-r$, the $\ell,r+t$-entry of $G$ is
$\alpha_{\ell t}$ where $\alpha_{\ell t}$ is 1 if $\ell=a_t$, 
$\alpha_{\ell t}=0$ if $\ell< a_t$ or $\ell=a_i$ (some $i\not=t$), and otherwise
$\alpha_{\ell t}\in\{0,1,\ldots,p-1\}$.
(So $G\in GL_n(\Z)$ with $G^{-1}=E_1E_2$
where $E_1$ has $i,i$-entry 1 ($1\le i\le n$), $\ell,r+t$-entry
$-\alpha_{d_{\ell}t}$, with zeros elsewhere, and $E_2$ has
$t,d_t$-entry 1 for $1\le t\le r$, $r+t,a_t$-entry 1 for $1\le t\le
n-r$, and zeros elsewhere.)  Let $\Y_r$ be the collection of matrices
$\pmatrix Y\\&0\endpmatrix$ where $Y$ varies over integral $r\times
r$, symmetric matrices modulo $p$, and let $D_r=\pmatrix
I_r\\&pI_{n-r}\endpmatrix$.

\proclaim{Theorem 2}  Given a degree $n$ Siegel modular form $F$,
$$F|T(p)=p^{n(k-n-1)/2}\sum_M F|\pmatrix {1\over p}I_n\\d&I_n\endpmatrix M$$
where $$M=\pmatrix D&Y\\&D^{-1}\endpmatrix \pmatrix
G^{-1}\\&^tG\endpmatrix$$
varies so that for some $r$, $0\le r \le n$, $D=D_r$, $Y\in\Y_r$, and
$G\in\G_r$.
\endproclaim

\demo{Proof}  Using Proposition 3.1 of [3], we only need to show that
as $G$ varies over $\G_r$, $\O=\L GpD_r^{-1}$ varies once over all lattices
$\O$ where $p\L\subseteq\O\subseteq\L$, $[\L:\O]=p^{r}$.  So, similar
to the proof of Theorem 1, we construct all the $\O$ as well as a
specific basis for each $\O$.

Let $\overline C$ be a dimension $n-r$ subspace of $\L/p\L$.  Choose a
basis $\overline v_1,\ldots,\overline v_{n-r}$ so that
$$\overline v_t=\overline x_{a_t}+\sum_{\ell>a_t}
\overline\alpha_{\ell t}\overline x_{\ell}$$
where $a_1,\ldots,a_{n-r}$ are distinct, $\overline\alpha_{\ell t}=0$
if $\ell=a_i$ (some $i\not=t$); for each $\overline\alpha_{\ell t}$,
 take a preimage
$\alpha_{\ell t}\in\{0,1,\ldots, p-1\}$.
Then with $(\{d_1,\ldots,d_r\},\{a_1,\ldots,a+{n-r}\})$ an ordered
partition
of $\{1,2,\ldots,n\}$
and $G$ constructed according to our recipe preceding
Theorem 2, we have $\O=\L GpD_r^{-1}$. $\square$
\enddemo

\remark{Remark}  In [6] we discuss how a particular subgroup of
$Sp_n(\Z)$ acts on Jacobi forms on
$f:\H_{(n-m)}\times\C^{n-m,m}\to\C$.  From this we see that for $1\le
j\le n-m$, 
$$f|T_j(p^2)=
\sum_M f|\pmatrix {1\over
p}I_j\\&I_{n-j}\\&&pI_j\\&&&I_{n-j}\endpmatrix M$$
where $$M=\pmatrix D&Y\\&D^{-1}\endpmatrix \pmatrix G^{-1}\\&^tG\endpmatrix$$
varies so that for some $r_0,r_2$ with $r_0+r_2\le j$,
$D=D_{r_0,r_2}$, $Y\in\Y'_{r_0,r_2}$, and $G\in\G_{\P}$ where $\P$ is
a partition of type $(r_0,r_2)$ so that
$\{n-m+1,\ldots,n\}\subseteq\{c_1,\ldots,c_{n-j}\}$. 
When $n-m<j\le n$, the operator $T_j(p^2)$ changes the index, as does $T(p)$, 
 and the matrices
giving the action of these operators is a bit more complicated;
explicit matrices for these operators are given in [6].
\endremark

\enddocument